\documentclass[10pt, letterpaper]{article}

\usepackage[dvips]{epsfig}
\usepackage{amssymb}
\usepackage{array}

\def\Fix{{\rm Fix}}

\newcommand{\binom}[2]{\left({{#1} \atop {#2}}\right)}

\newtheorem{Thm}{Theorem}
\newtheorem{Prop}[Thm]{Proposition}

\newtheorem{Cor}[Thm]{Corollary}

\newcommand{\desc}[2]{\begin{description}\item[{#1}]\ \ {#2} \end{description}\vspace{1mm}}

\begin{document}

\title{Enumeration of Symmetry Classes of Parallelogram
  Polyominoes\thanks{Work partially supported by NSERC (Canada) and
    FCAR (Qu\'ebec).}} 
\author{\textsc{P. Leroux, E. Rassart}}
\maketitle
\centerline{\textsl{LaCIM, Department of Mathematics}}
\centerline{\textsl{Universit\'e du Qu\'ebec \`a Montr\'eal} (UQAM)}

\begin{abstract}
Parallelogram polyominoes are a subclass of convex polyominoes in the
square lattice that has been studied extensively in the
literature. Recently \cite{PLERAR} congruence classes of convex
polyominoes with respect to
rotations and reflections have been enumerated by counting orbits
under the action of the dihedral group $\mathfrak{D}_4$, of symmetries of the
square, on (translation-type) convex polyominoes. Asymmetric
convex polyominoes were also enumerated using M\"obius inversion
in the lattice of subgroups of $\mathfrak{D}_4$. In this paper
we extend these results to the subclass of parallelogram polyominos
using a subgroup $\mathfrak{D}_2$ of $\mathfrak{D}_4$ which acts of this class.
\end{abstract}

\renewcommand\abstractname{\protect R\'esum\'e}

\begin{abstract}
Les polyominos parall\'elogrammes forment une sous-classe des
polyominos convexes sur le r\'eseau carr\'e; cette classe a \'et\'e
\'etudi\'ee en d\'etail dans la litt\'erature. Dans un travail
r\'ecent \cite{PLERAR}, les classes de congruences des polyominos
convexes sous les rotations et les r\'eflexions ont \'et\'e \'enum\'er\'ees
en d\'enombrant les orbites de l'action du groupe di\'edral
$\mathfrak{D}_4$, le groupe des sym\'etries du carr\'e, sur les
polyominos convexes (\`a translation pr\`es). L'inversion de M\"obius
dans le treillis des sous-groupes de $\mathfrak{D}_4$ a aussi
permis l'\'enum\'eration des polyominos convexes asym\'etriques. Nous
nous proposons dans cet article d'\'etendre ces r\'esultats \`a la
sous-classe des polyominos parall\'elogrammes, en faisant agir sur
elle un sous-groupe ($\mathfrak{D}_2$) de $\mathfrak{D}_4$. 
\end{abstract}

\section{Introduction}

\begin{figure}[!ht]
\begin{center}
\includegraphics{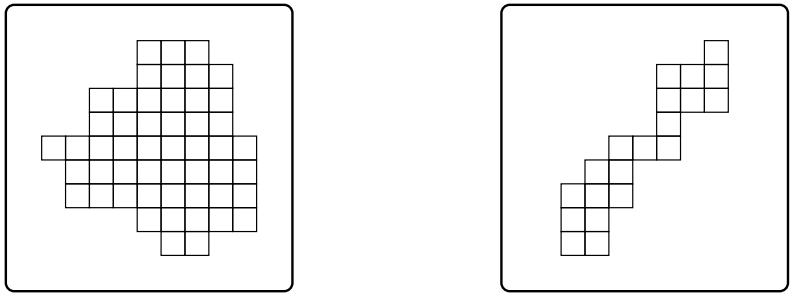}
\caption{\small Convex (left) and parallelogram (right) polyominoes.}
\label{ConvParfig}
\end{center}
\end{figure}

Parallelogram polyominoes, sometimes called staircase polyominoes,
form a subclass of (horizontally and vertically) convex polyominoes on
the square lattice, characterized by the fact that they touch the
bottom-left and the top-right corners of their minimal bounding
rectangle. See Figure~\ref{ConvParfig}. A $90^\circ$ rotation of these
would give a distinct but equivalent class of parallelogram
polyominoes. In the same way as for general convex polyominoes, the
area of a parallelogram polyomino is defined as the number of cells
that it contains and the half-perimeter is equal to the sum of its
width and height. Considerable literature can be found on the
enumeration of various classes of polyominoes having some convexity
and directedness property, with motivation coming from combinatorics,
statistical physics, computer science and recreational
mathematics. See M.~Bousquet-M\'elou \cite{MBMc} for a recent survey.
In particular, parallelogram polyominoes have been studied with
respect to their perimeter and area first by P\'olya, with further
contributions yielding refined enumerations from Bender, Klarner, Rivest,
Delest, F\'edou, Viennot and others. Mireille
Bousquet-M\'elou, using the Temperley methodology (\cite{HNVT}), has
given a generating function with respect to height, width, area and
height of first and last columns (\cite{MBMa}, \cite{MBMc}).

Polyominoes are usually considered equivalent if they can be obtained
from one another by a plane translation. They are sometimes called
translation-type polyominoes to be more precise (see D.~A.~Klarner
\cite{DAK}). It is natural to consider also \textsl{congruence-type}
polyominoes, that is, equivalence classes of polyominoes under
rotations and reflections. They occur as pieces that can
freely move in space, as in plane packing problems (see S.~W.~Golomb
\cite{SWG}). In \cite{PLERAR}, the enumeration of congruence-type
polyominoes, according to area and perimeter, has been carried out in
the case of convex polyominoes. 

The problem is equivalent to the
enumeration of orbits of the dihedral group $\mathfrak{D}_4$, of
symmetries of the square, acting on convex polyominoes. The group
$\mathfrak{D}_4$ contains eight elements, usually represented as $1$,
$r$, $r^2$, $r^3$, $h$, $v$, $d_1$ and $d_2$, where $1$ denotes the
identity element, $r$ denotes a rotation by a right angle, $h$ and $v$,
reflections with respect to the horizontal and vertical axes
respectively, and $d_1$ and $d_2$, reflections about the two diagonal
axes of the square (we take the bissector of the first quadrant for
$d_2$). The number of orbits $|X/G|$ of any finite group $G$ acting on
a set $X$ is given by the Cauchy-Frobenius formula (alias Burnside's
Lemma):
\begin{equation}
|X/G| = \frac{1}{|G|}\sum_{g \in G}|\mathrm{Fix}(g)|\,,
\label{eqnBurnside}
\end{equation}
where $\mathrm{Fix}(g)$ denotes the set of elements of X which are
$g$-symmetric, that is, invariant under $g$. Hence the enumeration of
congruence-type convex polyominoes involves determining the size of
the symmetry classes of convex polyominoes for each group element $g
\in \mathfrak{D}_4$. Formula~(\ref{eqnBurnside}) is valid for infinite
sets provided a weighted cardinality $|X|_{\omega}$ is taken, with
respect to some $G$-invariant weight function $\omega$. For a class
$\mathcal{P}$ of polyominoes this means using generating series
$\mathcal{P}(t,q)$ with respect to half-perimeter and area (variables
$t$ and $q$), for example.

The main goal of this paper is to carry out a similar procedure for
the class $\mathbb{P}$ of parallelogram polyominoes. We observe that a
subgroup of $\mathfrak{D}_4$ acts on parallelogram polyominoes, which
we denote $\mathfrak{D}_2$, namely $\mathfrak{D}_2 = \langle
r^2,d_1\rangle = \{1,r^2,d_1,d_2\}$, and that congruence types of
parallelogram polyominoes coincide with orbits of $\mathbb{P}$ under
$\mathfrak{D}_2$. In the following sections we therefore compute the
generating series of the symmetry classes $\mathrm{Fix}(g)$ of
parallelogram polyominoes for all $g \in \mathfrak{D}_2$ except the identity. We then use
(\ref{eqnBurnside}) to obtain $(\mathbb{P}/\mathfrak{D}_2)(t,q)$.

It is also possible to count asymmetric parallelogram polyominoes, that is
polyominoes that are not $g$-invariant for any $g$ except the
identity, using M\"obius
inversion in the lattice of subgroups of $\mathfrak{D}_2$. This
requires also the enumeration of the subclass
$\mathrm{Fix}(\mathfrak{D}_2)$ of $\mathbb{P}$, of totally symmetric
parallelogram polyominoes. We carry out this computation and show that
asymmetric parallelogram polyominoes are asymptotically equivalent to
all parallelogram polyominoes, as expected.

As we will see, the enumeration of all the symmetry classes of
parallelogram polyominoes, according to perimeter, involves in one way
or the other either the Dyck paths (or Dyck words, see J.~Labelle \cite{JLab}), counted by the Catalan numbers $c_n$, or the left
factors of Dyck paths, counted by the central binomial coefficients
$b_n$ (see Cori and Viennot \cite{RCXGV}), where
\begin{equation}
b_n = \binom{2n}{n} \qquad \textrm{and} \qquad c_n = \frac{1}{n+1}\binom{2n}{n}\,.
\end{equation}

When the area is taken into account, $q$-analogues (some well-known
and some novel) of these numbers appear naturally.

We would like to thank X.~G.~Viennot and M.~Bousquet-M\'elou for
useful discussions.

\section{Enumeration of parallelogram polyominoes}

It has been known for a long time (Levine \cite{JL}, P\'olya \cite{GP})
that the number of parallelogram polyominoes of perimeter $2n$ is
given by the Catalan number $c_{n-1} =
\frac{1}{n}\binom{2n-2}{n-1}$. One proof of this fact is provided by
the following bijection, due to Delest and Viennot (\cite{MPDXGV})
between parallelogram polyominoes of perimeter $2n+2$ and Dyck paths
of length $2n$: given a parallelogram polyomino $P$ of perimeter
$2n+2$, let 
$(a_1,a_2,\ldots,a_k)$ be the sequence of column heights of $P$, and
$(b_1,b_2,\ldots,b_{k-1})$ be such that $b_i$ is the number of cells
of contact between columns $i$ and $i+1$ of $P$. The associated Dyck
path $D$ is the unique Dyck path
with $k$ peaks and $k-1$ valleys such that the peak heights are given
in order by the sequence $(a_1,a_2,\ldots,a_k)$, and the valley heights
by the sequence $(b_1-1,b_2-1,\ldots,b_{k-1}-1)$ (the horizontal axis
is at level $0$). The height of $P$ is $n+1-k = n-(k-1)$, which is also given by
\begin{displaymath}
a_1 + (a_2-b_1) + (a_3-b_2) + \ldots + (a_k-b_{k-1}) =
\sum_{i=1}^ka_i - \sum_{j=1}^{k-1}b_j\,.
\end{displaymath}
On the other hand, the number of $\searrow$ steps in $D$, that is the half-length of
the path, is given by
\begin{displaymath}
\sum_{i=1}^ka_i - \sum_{j=1}^{k-1}(b_j-1) = \sum_{i=1}^ka_i -
\sum_{j=1}^{k-1}b_j + (k-1)\,,
\end{displaymath}
which is seen to be $n$ by the previous equation. Hence $D$ is a Dyck
path of length $2n$.

\begin{figure}[!ht]
\begin{center}
\includegraphics{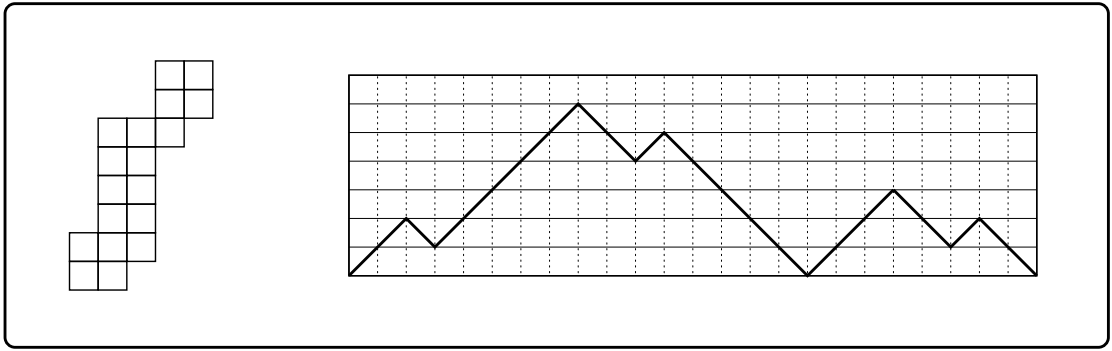}
\caption{\small Parallelogram polyomino and associated Dyck path.}
\label{DyckParBij1fig}
\end{center}
\end{figure}

Also, note that the sum of the heights of the peaks,
$\sum_{i=1}^ka_i$ is simply the area of
$P$. Figure~\ref{DyckParBij1fig} illustrates the 
bijection. 

It follows that the generating series $\mathbb{P}(t)$ for
parallelogram polyominoes according to half-perimeter is
\begin{equation}
\mathbb{P}(t) = \sum_{n \geq 2}c_{n-1}t^n = \frac{1-2t-\sqrt{1-4t}}{2}\,.
\end{equation}

It has also been known for some time that when the area is taken into
account, the generating series involves a quotient of $q$-analogues of
Bessel functions (see Klarner and Rivest \cite{DAKRLR} and Bender
\cite{EAB}). P\'olya (\cite{GP} and \cite{PFa}) found a Laurent series
relating the area and perimeter generating function to a
specialization of itself, from which the terms of the series can be
extracted easily. The width, height and area generating series can
also be expressed as a continued fraction (see \cite{PFb}). We will
use the following recent more general form due to
M.~Bousquet-M\'elou~\cite{MBMa}, giving the generating series
$\mathbb{P}(v,x,y,q)$ of parallelogram polyominoes, where the
variables $v$, $x$, $y$ and $q$ mark respectively the height of the
rightmost column, the width, the overall height, and the area.

\begin{Prop}\textnormal{(\cite{MBMa})}
The generating function $\mathbb{P}(v,x,y,q)$ of parallelogram
polyominoes is given by
\begin{equation}
\mathbb{P}(1,v,x,y,q) = vy\frac{J_1(v,x,y,q)}{J_0(x,y,q)}\,,
\label{eqnMBMPar}
\end{equation}
with
\begin{equation}
J_0(x,y,q) = \sum_{n \geq 0}\frac{(-1)^nx^nq^{\binom{n+1}{2}}}{(q)_n(yq)_n} 
\end{equation}
and
\begin{equation}
J_1(v,x,y,q) = \sum_{n \geq
  1}\frac{(-1)^{n-1}x^nq^{\binom{n+1}{2}}}{(q)_{n-1}(yq)_{n-1}(1-vyq^n)} 
\end{equation}
with the usual notation $(a)_n = (a;q)_n = \prod_{i=0}^{n-1}(1-aq^i)$.
\hfill$\blacksquare$
\end{Prop}

Note that the half-perimeter and area generating function
$\mathbb{P}(t,q)$ of parallelogram polyominoes is obtained by putting
$v=1$, $x=t$, $y=t$ in (\ref{eqnMBMPar}).

\section{Symmetry classes of parallelogram polyominoes}

\subsection{Rotational symmetry}

Observe that if we apply the Delest-Viennot bijection to an
$r^2$-symmetric parallelogram polyomino of perimeter $2k+2$, the Dyck
path we obtain is vertically symmetric (or, equivalently, the Dyck
word associated to it is a palindrome). Hence we need only consider
half the path, which is simply a left factor, of length $k$, of a Dyck
path. 

\begin{Prop}
The number of $r^2$-symmetric parallelogram polyominoes of
half-perimeter $k+1$ is equal to the number of left factors of Dyck
paths, of length $k$.\hfill$\blacksquare$
\end{Prop}

\begin{Cor}The number of $r^2$-symmetric parallelogram polyominoes of
  half-perimeter $k+1$ is given by 
\begin{equation}
r_{k+1}(1) = \left\{\begin{array}{ll}
\binom{k}{k/2} & \textrm{if $k$ is even}\,,\\[2mm]
\frac{1}{2}\binom{k+1}{(k+1)/2} & \textrm{if $k$ is odd}\,.
\end{array}\right.
\label{eqn:r2symnum}
\end{equation}
\label{cor:r2symnum}
\end{Cor}

\desc{Proof}{It is known that the number of left factors of length
  $2n$ of Dyck words is equal to the number of words in the alphabet \{0,1\} with distribution
  $0^n1^n$, from which (\ref{eqn:r2symnum}) follows
  easily. See \cite{RCXGV} for a bijective proof. Here we prove~(\ref{eqn:r2symnum}) using generating functions. Dyck paths and left factors of
  Dyck paths are generated by the algebraic grammar 
\begin{eqnarray*}
C \rightarrow \varepsilon + xC\bar{x}C\\
L \rightarrow C + CxL\,,
\end{eqnarray*}
over the alphabet $\{x, \bar{x}\}$. $C$ denotes the Dyck paths and $L$
the left factors, while $x$ and $\bar{x}$ respectively denote a
$\nearrow$ step and a $\searrow$ step. The first production rule gives
$C(x,\bar{x}) = 1 + x\bar{x}C(x,\bar{x})^2$, which we solve for 
\begin{displaymath}
C(x,\bar{x}) = \frac{1-\sqrt{1-4x\bar{x}}}{2x\bar{x}}\,.
\end{displaymath} 

The second production rule gives $L(x,\bar{x}) = C(x,\bar{x})(1 +
xL(x,\bar{x}))$, which we can solve, now that we have $C(x,\bar{x})$,
for 
\begin{displaymath}
L(x,\bar{x}) = \frac{\sqrt{1-4x\bar{x}}-1}{x(1-2\bar{x}-\sqrt{1-4x\bar{x}})}\,.
\end{displaymath}

Substituting $x \mapsto t, \bar{x} \mapsto t$ into $L(x,\bar{x})$
gives the generating series $L(t)$ of left factors of Dyck paths by
their length: 
\begin{displaymath}
L(t) = \frac{2t-1+\sqrt{1-4t^2}}{2t(1-t)}\,,
\end{displaymath}
from which~(\ref{eqn:r2symnum}) follows.
\hfill$\blacksquare$
}

\begin{figure}[!ht]
\begin{center}
\includegraphics{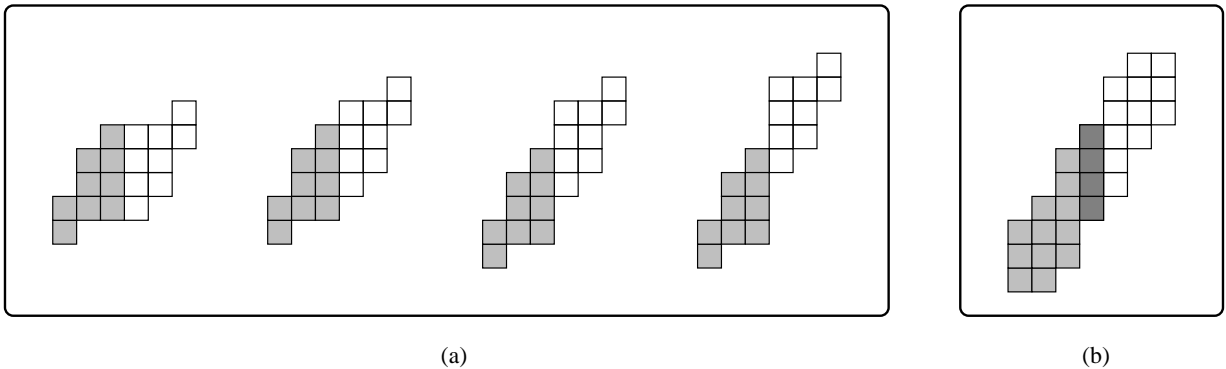}
\caption{\small $r^2$-symmetric parallelogram polyominoes with (a) even width and (b) odd width.}
\label{r2symfig}
\end{center}
\end{figure}

In order to include the area, we could extend another bijection, due
to Bousquet-M\'elou and Viennot \cite{MBMXGV}, involving heaps of
segments, to left factors of Dyck paths. However, there is a more direct approach. Indeed,
parallelogram polyominoes with rotational symmetry can be obtained
from two copies of a same parallelogram polyomino glued together. The
glueing process depends on whether we want the final object to be of
even width or of odd width, as can be seen in
Figure~\ref{r2symfig}. If $R_2(x,y,q)$ is the generating function of
$r^2$-symmetric parallelogram polyominoes, then
\begin{equation}
R_2(x,y,q) = R_2^{(e)}(x,y,q) + R_2^{(o)}(x,y,q)
\end{equation}
where $R_2^{(e)}(x,y,q)$ and $R_2^{(o)}(x,y,q)$ are respectively the
generating series of even-width and odd-width $r^2$-symmetric
parallelogram polyominoes.  

\begin{Prop}
The generating function $R_2^{(e)}(x,y,q)$ of even-width parallelogram polyominoes is given by
\begin{equation}
\label{eqn:r2sym}
R_2^{(e)}(x,y,q) = \frac{1}{1-y} \left(\mathbb{P}(\frac{1}{y}, x^2,
  y^2, q^2) - \mathbb{P}(1, x^2, y^2, q^2) \right)\,. 
\end{equation}
where $\mathbb{P}(v,x,y,q)$ is the generating function~\textnormal{(\ref{eqnMBMPar})} of parallelogram polyominoes.
\end{Prop}

\desc{Proof}{Let $P$ be an $r^2$-symmetric parallelogram of even
  width. We define the fundamental region of $P$ to be the left half
  of $P$. Call this polyomino $Q$ (see Figure~\ref{r2symfig}(a)). We
  first remark that $Q$ is a parallelogram polyomino. To get $P$ from
  $Q$, we rotate a copy of $Q$ by $180^\circ$ and glue the result
  $\overline{Q}$ to $Q$ along the rightmost column. If this column has
  length equal to $k$, there will be $k$ possible positions for
  $\overline{Q}$ relative to $Q$. The substitution $v \mapsto 1/y$, $x
  \mapsto x^2$, $y \mapsto y^2$ and $q \mapsto q^2$ in the generating
  series $P(1,v,x,y,q)$ of directed convex polyominoes corresponds to
  the highest position of $\overline{Q}$, which minimizes the overall
  height of $P$. All the possible positions will be accounted for by
  multiplying by $(1+y+\ldots + y^{k-1})$. In other words, the
  substitution to make in $\mathbb{P}(v,x^2,y^2,q^2)$ is 
\begin{equation}
v^k \mapsto \frac{1 + y + \ldots + y^{k-1}}{y^k} = \frac{1}{1-y}\left(\frac{1}{y^k}-1\right)\,.
\end{equation}

Summing over all possible $k$'s, we find the proposed generating
series~(\ref{eqn:r2sym}) for $r^2$-symmetric parallelograms. 
\hfill$\blacksquare$ 
}

\begin{Prop}
The generating function $R_2^{(o)}(x,y,q)$ of odd-width parallelogram polyominoes is given by
\begin{equation}
R_2^{(o)}(x,y,q) = \frac{1}{x}\mathbb{P}(\frac{1}{yq}, x^2, y^2, q^2)\,.
\end{equation}
where $\mathbb{P}(v,x,y,q)$ is the generating function of parallelogram polyominoes.
\end{Prop}

\desc{Proof}{The proof is similar to the previous one. The main
  difference is that only one glueing position of $\overline{Q}$ to
  $Q$ is admissible and that furthermore the rightmost column of $Q$
  and its rotated image in $\overline{Q}$ are superimposed to yield an
  odd width (see Figure~\ref{r2symfig}(b)). Details are left to the
  reader. 
\hfill$\blacksquare$
}

We would like to find the number of $r^2$-symmetric parallelograms of
a given half-perimeter, without losing the area information, i.e. we
want to express the generating series in the form 
\begin{equation}
R_2(t,q) = R_2(t,t,q) = \sum_{k \geq 0}r_k(q)t^k\,. 
\label{eqnr2kq}
\end{equation}

The above expressions for the generating series of $r^2$-symmetric
parallelogram polyominoes can be used to extract the 
polynomials $r_k(q)$ from it (i.e. develop it in powers of $t$ after
substituting $x \mapsto t, y \mapsto t$ in it). Here
are the first few of these polynomials: 
\begin{eqnarray*}
r_2(q) & = & q\\
r_3(q) & = & 2q^2\\
r_4(q) & = & q^4 + 2q^3\\
r_5(q) & = & 2q^6 + 2q^4\\
r_6(q) & = & q^9 + 2q^8 + q^7 + 2q^6 + 4q^5
\end{eqnarray*}

\subsection{Reflective symmetries}

\begin{figure}[!ht]
\begin{center}
\includegraphics{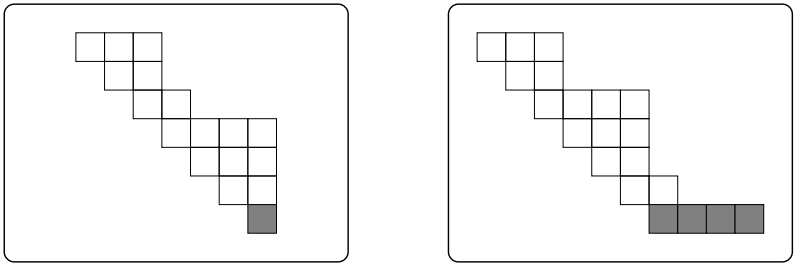}
\caption{\small Dyck polyomino (left) and left factor of Dyck polyomino (right).}
\label{DyckPolyofig}
\end{center}
\end{figure}

We begin by introducing a subfamily of parallelogram polyominoes which
we will call \textsl{Dyck polyominoes} as they correspond to Dyck
paths drawn over and above the main diagonal. We will also consider
truncated Dyck polyominoes, which we will
call \textsl{left factors of Dyck polyominoes} (or $LFD$ polyominoes
for short), again in analogy with the left factors of Dyck paths. Dyck
and $LFD$ polyominoes are illustrated in Figure~\ref{DyckPolyofig}.

\begin{figure}[!ht]
\begin{center}
\includegraphics{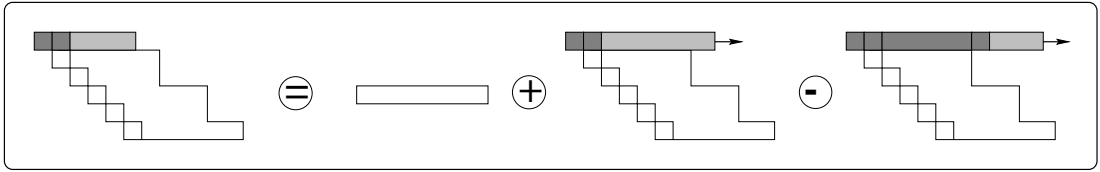}
\caption{\small Recursive construction of the $LFD$ polyominoes, by adding rows..}
\label{LFDconstructfig}
\end{center}
\end{figure}

We introduce $L_n(u) = L_n(u,y,q)$ the generating function for $LFD$
polyominoes having a basis of width $n$, with variables $u$, $y$ and
$q$ corresponding to the number of cells of the uppermost row,
the height and the area respectively. $L_n(u)$ can be defined recursively by the
following functional equation, illustrated in Figure~\ref{LFDconstructfig}:
\begin{equation}
L_n(u) = u^nyq^n + \frac{yu^2q^2}{1-uq}(L_n(1)-L_n(uq))\,.
\label{eqnLFD}
\end{equation}

The generating function $L(u)$ of all $LFD$ polyominoes is simply the sum over all possible base widths,
\begin{equation}
L(u) = \sum_{n \geq 1}L_n(u)\,.
\end{equation}

Moreover, the Dyck polyominoes being the $LFD$ polyominoes with width
one bases, their height and area generating series $D(t,q)$ is given by
\begin{equation}
D(y,q) = L_1(1,y,q)\,.
\end{equation}

A straightforward application of Lemma 2.3 from \cite{MBMa} (M. Bousquet-M\'elou) gives the solution to the functional equation~(\ref{eqnLFD}). As we do not need the variable $u$ for our purpose, we set it equal to $1$, which simplifies the expression for the generating function.

\begin{Prop}
The area and height generating function $L_n(1,y,q)$ for $LFD$
polyominoes having a basis of width $n$ is given by
\begin{equation}
L_n(1,y,q) = \frac{\displaystyle\sum_{m \geq
    0}\frac{(-1)^my^{m+1}q^{(m+n)(m+1)}}{(q)_m}}{\displaystyle
  \sum_{m \geq 0}\frac{(-1)^my^{m}q^{m(m+1)}}{(q)_{m}}}\,.
\label{eqnLn}
\end{equation}\hfill$\blacksquare$
\end{Prop}

For $n=1$, this gives 
\begin{equation}
D(y,q) = \frac{\displaystyle\sum_{m \geq
    0}\frac{(-1)^my^{m+1}q^{(m+1)^2}}{(q)_m}}{\displaystyle
  \sum_{m \geq 0}\frac{(-1)^my^{m}q^{m(m+1)}}{(q)_{m}}}
\end{equation}
for the height and area generating function for Dyck
polyominoes. However, we can also express $D(y,q)$ using the
classical $q$-analogue of Catalan numbers $c_n(q)$, satisfying
the recurrence 
\begin{equation}
c_n(q) = \sum_{k=0}^{n-1}q^kc_k(q)c_{n-1-k}(q)\,,
\end{equation}
as it is well known that they area-enumerate Dyck paths of length
$2n$. The area enumerated by $c_n(q)$ is the number of cells under the
path and strictly above its supporting diagonal (i.e. the cells on the
diagonal are not included in the area). To get a Dyck polyomino from a
Dyck path, we have to add the area of the diagonal. If the length of
the path is $2n$, then a factor of $q^n$ has to be added. A further
diagonal of cells has to be added because the Dyck paths can touch the
supporting diagonal, in which case they are not polyominoes. For the
paths of length $2n$, $n+1$ cells thus have to be added, contributing a
further $q^{n+1}$ factor to the area. This last diagonal
also adds one unit of height to the polyominoes. Hence
\begin{equation}
D(y,q) = \sum_{n \geq 1}y^nq^{2n-1}c_{n-1}(q)\,.
\end{equation}

\subsubsection{Reflective symmetry along the first diagonal}

\begin{figure}[!ht]
\begin{center}
\includegraphics{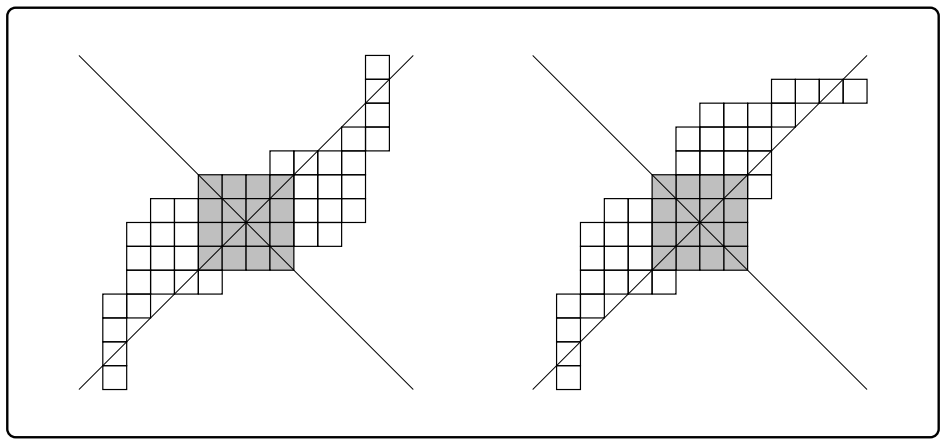}
\caption{\small Bijection between $r^2$-symmetric parallelogram
  polyominoes of even half-perimeter and $d_1$-symmetric
  parallelograms.} 
\label{r2d2bijfig}
\end{center}
\end{figure}

There is a nice area-preserving
bijection between $d_1$-symmetric parallelograms of a given
half-perimeter and $r^2$-symmetric parallelograms with same
half-perimeter. Since the minimal rectangle of a $d_1$-symmetric
parallelogram is necessarily a square, the perimeter is a multiple of $4$, and thus the
half-perimeter is even. Hence we have 
\begin{equation}
D_1(x,y,q) = \sum_{k \geq 0}r_{2k}(q)t^{2k}\,,
\end{equation}
where $D_1(x,y,q)$ is the generating series of $d_1$-symmetric
parallelogram polyominoes and the $r_{2k}(q)$ are defined by~(\ref{eqnr2kq}). The bijection is shown on an example in
Figure~\ref{r2d2bijfig}, and goes as follows: a $r^2$-symmetric
parallelogram has a center of rotation. If it has even half-perimeter,
this center will either fall in the center of a cell (if both the
height and the width are odd) or be the common corner of four cells
forming a square (if both the height and the width are even). In both
cases, we consider the first diagonal (parallel to the bissector of
the second quadrant) passing through the center of rotation and the
region of the parallelogram below the second diagonal (see
Figure~\ref{r2d2bijfig}). This region is not a polyomino, but the
parallelogram is obtained by glueing the region and a copy of it
rotated by $180^\circ$ in the unique way such that there are no
``half-cells'' left. Suppose that instead of rotating the copy of the
region, we reflect it along the second diagonal and glue it so that
there are no half-cells left, then we clearly obtain a $d_1$-symmetric
parallelogram which, further, has the exact same perimeter and area as
the initial $r^2$-symmetric parallelogram. We can similarly reverse
the process to start with an arbitrary $d_1$-symmetric parallelogram
and end with a $r^2$-symmetric parallelogram. 

\subsubsection{Reflective symmetry along the second diagonal}

\begin{figure}[!ht]
\begin{center}
\includegraphics{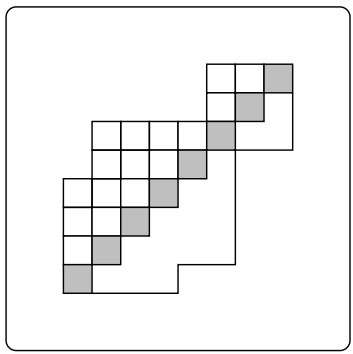}
\caption{\small A $d_2$-symmetric parallelogram polyomino.}
\label{d1symfig}
\end{center}
\end{figure}

We next consider $d_2$-symmetric parallelogram polyominoes,
i.e. parallelograms which are left invariant by a symmetry along the
second diagonal. Figure~\ref{d1symfig} gives an example of such a
parallelogram. We observe first that the minimal rectangle of such a
parallelogram will always be a square with side length equal to the
quarter-perimeter of the inscribed parallelogram.

We note that $d_2$-symmetric parallelogram polyominoes can be
constructed from two copies of a same Dyck polyomino, whose diagonals
we glue together (dark cells on Figure~\ref{d1symfig}). The area of
the final object will be twice the area of the Dyck polyomino minus
the area of diagonal, which was counted twice. There are as many cells
on the diagonal as the height of the Dyck polyomino, and the width of
the final object will also be the height of the Dyck polyomino. Hence we get

\begin{Prop}
The generating series $D_2(x,y,q)$ of $d_2$-symmetric parallelogram polyominoes is given by
\begin{equation}
D_2(x,y,q) = D(\frac{xy}{q},q^2)\,.
\end{equation}
\hfill$\blacksquare$
\end{Prop}

\subsubsection{Reflective symmetry along both diagonals}

\begin{figure}[!ht]
\begin{center}
\includegraphics{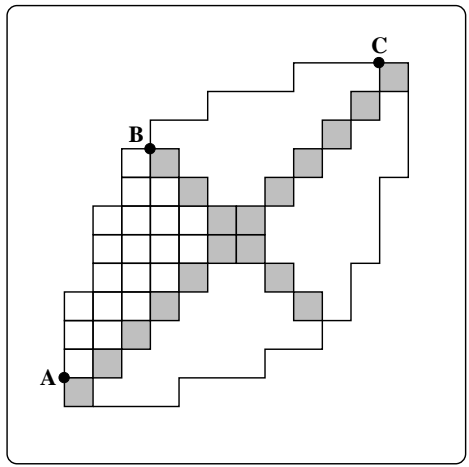}
\caption{\small A $\mathfrak{D}_2$-symmetric parallelogram polyomino.}
\label{r2d1d2symfig}
\end{center}
\end{figure}

The final (non-cyclic) subgroup whose set of fixed elements we study is the
whole group itself. This group is generated by any two nontrivial
elements, but it is convenient to consider the symmetries along the
two diagonals as the generators. This allows us to characterize the
fundamental region of a $\mathfrak{D}_2$-symmetric parallelogram, as
can be seen in Figure~\ref{r2d1d2symfig}. 

We note first that the minimal rectangle of a
$\mathfrak{D}_2$-symmetric parallelogram $P$ is a square. We remark
also that the exterior path going from $\mathbf{A}$ to $\mathbf{C}$ is
a Dyck path that has the additional property that it is symmetric with
respect to the second diagonal passing through the center of $P$
(i.e. the Dyck word in $x$ and $\bar{x}$ associated to the Dyck path
is a palindrome). Hence $P$ is completely determined by ``half'' a
Dyck path (the path going from $\mathbf{A}$ to $\mathbf{B}$). If $P$
has half-perimeter $2k$ (its half-perimeter is necessarily even since the
minimal rectangle is a square), then the path $\mathbf{A} \rightarrow
\mathbf{C}$ is a symmetrical Dyck path of length $2k-2$, and the path
$\mathbf{A} \rightarrow \mathbf{B}$ is simply a left factor of length
$k-1$ of a Dyck path. Thus we have the following result: 

\begin{Prop}
The number of $\mathfrak{D}_2$-symmetric parallelogram polyominoes of half-perimeter $2k+2$ is given by
\begin{equation}
d^{(1,2)}_{2k+2}(1) =  \left\{\begin{array}{ll}
\binom{k}{k/2} & \textrm{if $k$ is even}\,,\\[2mm] 
\frac{1}{2}\binom{k+1}{(k+1)/2} & \textrm{if $k$ is odd}\,.
\end{array}\right.
\end{equation}
\end{Prop}

\desc{Proof}{See Corollary~\ref{cor:r2symnum}.\hfill$\blacksquare$
} 

\begin{figure}[!ht]
\begin{center}
\includegraphics{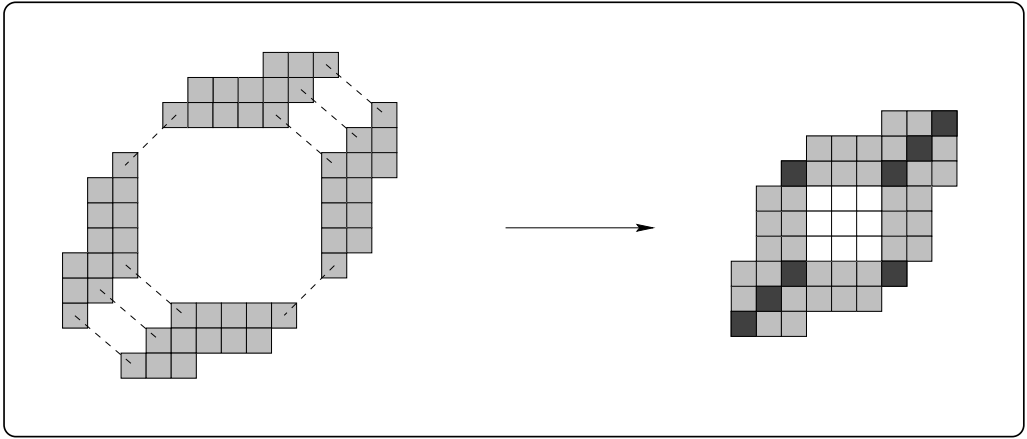}
\caption{\small Construction of a $\mathfrak{D}_2$-symmetric parallelogram polyomino from $4$ copies of a $LFD$ polyomino.}
\label{d1d2constructfig}
\end{center}
\end{figure}

We obtain the area and half-perimeter generating function for the
$\mathfrak{D}_2$-symmetric parallelogram polyominoes by constructing
these from 4 copies of a $LFD$ polyomino, as illustrated in
Figure~\ref{d1d2constructfig}. Some cells are superposed (the dark
ones) and others have to be added (square of white cells in the
center), and that has to be taken into account when computing the area
of the final object. If the $LFD$ polyomino has area $A$, height
(number of cells on the diagonal) $d$ and base $n$, then we have that
the area of the $\mathfrak{D}_2$-symmetric parallelogram polyomino is
$4A - 2d + (n-2)^2 - 2$, while its half-perimeter will be given by $2n
+ 4(d-1)$. Hence we have the following proposition:

\begin{Prop}
The half-perimeter and area generating function $D_{1,2}(t,q)$ for
$\mathfrak{D}_2$-symmetric parallelogram polyominoes is given by
\begin{equation}
D_{1,2}(t,q) = t^2q + \sum_{n \geq 2}t^{2n-4}q^{n^2-4n+2}L_n(1,\frac{t^4}{q^2},q^4)
\end{equation}
where $L_n(u,x,q)$ is the generating series for $LFD$ polyominoes with
a base of width $n$.\hfill$\blacksquare$ 
\end{Prop}

\begin{Cor}
\begin{equation}
D_{1,2}(t,q) = t^2q + \frac{\displaystyle \sum_{n \geq 2}\sum_{m \geq
  0}\frac{(-1)^mt^{4m+2n}q^{4m^2+2m+4mn+n^2}}{(q^4)_m}}{\displaystyle
  1-t^4q^6\sum_{m \geq 0}\frac{(-1)^mt^{4m}q^{4m^2+10m}}{(q^4)_{m+1}}}\,.
\end{equation}
\end{Cor}

\desc{Proof}{This follows from equation~(\ref{eqnLn}).\hfill$\blacksquare$}

Here are the first few terms of $D_{1,2}(t,q)$:
\begin{displaymath}
D_{1,2}(t,q) = t^2q + t^4q^4 + t^6q^9 + t^8q^{10} + t^8q^{14} + t^8q^{16} + t^{10}q^{15} + t^{10}q^{19} + t^{10}q^{23} + t^{10}q^{25} + \ldots
\end{displaymath}

\subsection{Congruence-type parallelogram polyominoes}

We are now in a position to enumerate congruence-type parallelogram
polyominoes, i.e. parallelograms up to rotation and reflection using
formula~(\ref{eqnBurnside}) with $G = \mathfrak{D}_2$ and $\mathcal{P}
= \mathbb{P}$, the class of all parallelogram polyominoes: 
\begin{equation}
|\mathbb{P}/\mathfrak{D}_2|_w = \frac{1}{4}\sum_{g \in \mathfrak{D}_2} |\Fix(g)|_w\,,
\end{equation}
where $|\Fix(g)|_w$ is the half-perimeter and area generating series
of the convex $g$-symmetric polyominoes. Therefore, 

\begin{Prop}The half-perimeter and area generating series
  $(\mathbb{P}/\mathfrak{D}_2)(t,q)$ of congruence-type parallelograms is given
  by 
\begin{equation}
(\mathbb{P}/\mathfrak{D}_2)(t,q) = |\mathbb{P}/\mathfrak{D}_2|_w =
\frac{1}{4} \left(\mathbb{P}(1,t,t,q) + R_2(t,t,q) + D_1(t,t,q) +
  D_2(t,t,q)\right)\,. 
\end{equation}\hfill$\blacksquare$
\end{Prop}

Here are the first few terms of $(\mathbb{P}/\mathfrak{D}_2)(t,q) =
\sum_{k \geq 0}\widetilde{p}_k(q)t^k$:
\begin{eqnarray*}
\widetilde{p}_{2}(q) & = & q\\
\widetilde{p}_{3}(q) & = & q^{2}\\
\widetilde{p}_{4}(q) & = & q^{4} + 2\,q^{3}\\
\widetilde{p}_{5}(q) & = & q^{6} + q^{5} + 3\,q^{4}\\
\widetilde{p}_{6}(q) & = & q^{9} + 2\,q^{8} + 3\,q^{7} + 4\,q^{6} + 6\,q^{5}
\end{eqnarray*}

\subsection{Asymmetric parallelogram polyominoes}

We can also enumerate asymmetric parallelogram polyominoes,
i.e. parallelograms having no symmetry at all, using M\"obius
inversion in the lattice of subgroups of $\mathfrak{D}_2$.

\begin{figure}[!ht]
\begin{center}
\includegraphics{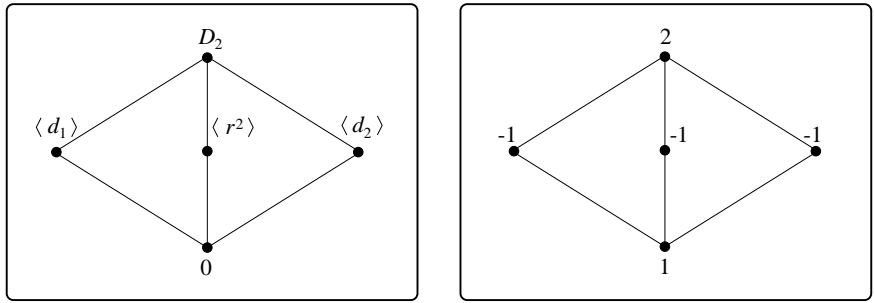}
\caption{\small Lattice of the subgroups of $\mathfrak{D}_2$ and the
  M\"obius function $\mu(0, H)$ on the lattice.}
\label{MobiusParfig}
\end{center}
\end{figure}

The reader is refered to \cite{GCR} for a general discussion of
M\"obius inversion, and to \cite{PLERAR} to see it applied
to the enumeration of the symmetry classes of convex polyominoes. We
simply give here in Figure~\ref{MobiusParfig} the lattice of the
subgroups of $\mathfrak{D}_2$ and
the value of the M\"obius function on the points of the lattice.
For subgroups $H$ of $\mathfrak{D}_2$, we denote by $F_{\geq H}$
(resp. $F_{= H}$) the half-perimeter and area generating series for the set of
parallelogram polyominoes having at least (resp. exactly) the
symmetries of $H$.

\begin{Prop}The half-perimeter and area generating series
  $\overline{\mathbb{P}}(t,q) = F_{=0}$ of asymmetric parallelogram
  polyominoes is given by 
\begin{eqnarray}
\overline{\mathbb{P}}(t,q) & = & F_{\geq 0} - F_{\geq \langle
  r^2\rangle} - F_{\geq \langle d_1\rangle} -
F_{\geq \langle d_2\rangle} + 2F_{\geq \langle d_1, d_2\rangle}\nonumber \\* 
& & \nonumber \\*
& = & \mathbb{P}(1,t,t,q) - R_2(t,t,q) - D_1(t,t,q) - D_2(t,t,q) + 2D_{1,2}(t,t,q)\,.\label{eqn:asym} 
\end{eqnarray}
where $D_{1,2}(x,y,q)$ is the generating series of $\mathfrak{D}_2$-symmetric polyominoes.
\hfill$\blacksquare$
\end{Prop}

Here are the first few terms of $\overline{\mathbb{P}}(t,q) = \sum_{k \geq 0}\overline{p}_k(q)t^k$:
\begin{eqnarray*}
\overline{p}_2(q) = \overline{p}_3(q) = \overline{p}_4(q) & = & 0\\
\overline{p}_5(q) & = & 4\,q^{5} + 4\,q^{4}\\
\overline{p}_6(q) & = & 8\,q^{7} + 8\,q^{6} + 8\,q^{5}\\
\overline{p}_7(q) & = & 4\,q^{11} + 8\,q^{10} + 20\,q^{9} + 24\,q^{8} + 32\,q^{7} + 24\,q^{6}
\end{eqnarray*}

Note that the same method would allow us to enumerate the convex polyominoes having exactly the symmetries of any given subgroup of $\mathfrak{D}_2$.

\subsection{Asymptotic results}

Here we show the asymptotic result that for large area or large perimeter,
almost all parallelogram polyominoes are asymmetric. In other words, the
probability for a parallelogram polyomino to have at least one symmetry goes
to zero as the area or the perimeter goes to infinity.

We know from the work of P\'olya (\cite{GP}) that the number of parallelogram
polyominoes of half-perimeter $n$ is given by
\begin{equation}
p_n^{(t)} = c_{n-1} \sim \frac{4^n}{\sqrt{\pi}n^{3/2}}\,.
\end{equation} 

For the area, we have the result from Bender:
\begin{Prop}\textnormal{(\textbf{Bender} \cite{EAB})}\ \ 
Let $p_n^{(q)}$ be the number of parallelogram polyominoes with area $n$. Then
\begin{equation}
p_n^{(q)} \sim k\,\mu^n\,, \label{eqn:asympt}
\end{equation}
with
\begin{displaymath}
k = 0.29745\ldots \qquad \mu = 2.30913859330\ldots 
\end{displaymath}\hfill$\blacksquare$
\end{Prop}

\begin{Prop}
Let $H$ be any non-trivial subgroup of $\mathfrak{D}_2$ and denote by
$P_{H}^{(q)}(n)$ (resp. $P_{H}^{(t)}(n)$) the number of $H$-symmetric
parallelogram polyominoes with area (resp. half-perimeter) $n$. Then, 
\begin{equation}
\lim_{n \rightarrow \infty}\frac{P_{H}^{(q)}(n)}{p_n^{(q)}} = 0 \qquad
\textrm{and} \qquad
\lim_{n \rightarrow \infty}\frac{P_{H}^{(t)}(n)}{p_n^{(t)}} = 0\,.
\end{equation}
\end{Prop}

\desc{Proof}{The proof for the perimeter part of the Proposition is immediate as we have
  closed forms for all the coefficients, and thus the limit can be
  verified to be zero explicitly. 

For the area, we need only consider $r^2$- and $d_2$-symmetric
parallelogram polyominoes, as the $d_1$-symmetric parallelograms are
in bijection with a subclass of $r^2$-symmetric ones, and the
$\mathfrak{D}_2$-symmetric parallelograms are a subclass of all the
other classes of symmetry. A same basic argument works for $r^2$- and
$d_2$-symmetric parallelogram polyominoes, using the fact that
they are constructed from two congruent subpolyominoes. A
supplementary column sometimes has to be added in the $r^2$ case, according
to whether the height or the width of the initial parallelogram is odd or
even. These subpolyominoes are parallelograms in every case, and they have at
most half the area of the initial object. 
\begin{itemize}
\item \textbf{$r^2$-symmetric parallelogram polyominoes\,:}\ An
$r^2$-symmetric parallelogram polyomino with even width and area $n$ is
constructed from two congruent subparallelograms with exactly half
the area, which can be glued together in at most $n/2$ ways (the
maximal height of the columns that get glued together). So
$P_{r^2}^{(q), even}(n) \leq \frac{1}{2}np_{n/2}^{(q)}$. Hence 
\begin{displaymath}
\lim_{n \rightarrow \infty}\frac{P_{r^2}^{(q),even}(n)}{p_n^{(q)}} \quad \leq
\quad  \lim_{n \rightarrow
\infty}\frac{\frac{1}{2}nk\mu^{n/2}}{k\mu^n} \quad = \quad 0\,. 
\end{displaymath}

Next consider an $r^2$-symmetric parallelogram polyomino with odd width and
area $n$. This polyomino is constructed from a central column ($n$
choices of height) and two congruent subparallelograms of area at
most $\lfloor n/2 \rfloor$. Then there are at most $n$ possible
positions where to glue the subparallelograms to the central column (they
are glued symmetrically). Thus $P_{r^2}^{(q), odd}(n) \leq
n^2(1+p_1^{(q)}+p_2^{(q)}+\ldots + p_{\lfloor n/2 \rfloor}^{(q)}) < n^3p_{\lfloor n/2
\rfloor}^{(q)}$ and the result follows as above. Hence the result holds for
the subgroup $\langle r^2 \rangle$ of $\mathfrak{D}_2$; 

\item \textbf{$d_2$-symmetric convex polyominoes\,:}\ Let
$P$ be a $d_2$-symmetric parallelogram polyomino and $Q$ its fundamental
region. Suppose that $P$ has $b$ cells on the diagonal symmetry
axis. Then the minimum area $P$ can have is $b + 2(b-1)$. This gives a
minimum area of $b + (b-1)$ for $Q$. Hence 
\begin{displaymath}
\frac{\textrm{Area of } Q_{min}}{\textrm{Area of } P_{min}} \quad =
\quad \frac{2b-1}{3b-2}\,.  
\end{displaymath}

Then if we add to $Q$ a cell not on the diagonal symmetry axis, two
cells get added to $P$, and thus we conclude that the ratio can only
decrease as we make $P$ into a larger $d_2$-symmetric parallelogram polyomino
with the same number of cells on the diagonal axis. For $b \geq 2$,
the ratio will be smaller than or equal to $3/4$. As a loose
approximation, we can take $Q$ to be any parallelogram polyomino. This gives
$P_{d_2}(n) \leq p_{\lceil 3n/4 \rceil}^{(q)} + 1$. The
$+1$ term corresponds to the unique $d_2$-symmetric parallelogram
polyominoes having only one cell on the diagonal axis. Hence
the result will also hold for the subgroup $\langle d_2 \rangle$.  
\end{itemize}\hfill$\blacksquare$
}

\begin{Prop}
If we denote by $\overline{p}_{n}^{(q)}$ (resp. $\overline{p}_n^{(t)}$)
the number of asymmetric parallelogram
polyominoes of area (resp. half-perimeter) $n$, then 
\begin{eqnarray}
\overline{p}_{n}^{(q)} & \sim & p_n^{(q)}\,,\\
\overline{p}_{n}^{(t)} & \sim & p_n^{(t)}\,.
\end{eqnarray}
\end{Prop}

\desc{Proof}{We get the result from equation~(\ref{eqn:asym}) and from
  the previous Proposition.\hfill$\blacksquare$ 
}

\vspace{5mm}

Two tables can be found in appendix that present the numbers of parallelogram polyominoes
according to their symmetry types and their perimeter or area. The
columns indexed by subgroups of $\mathfrak{D}_2$ give the numbers of
parallelogram polyominoes of a given perimeter or area that are left fixed by
the symmetries of the subgroup. The columns \textsl{\#\,Orbits} and
\textsl{Asym} give respectively the number of
congruence-type and asymmetric parallelogram polyominoes of the given size.

\pagebreak

\noindent{\LARGE Appendix}
\vspace{5mm}

\begin{table}[!ht]
\begin{center}
\begin{footnotesize}
\begin{tabular}{|c||r|r|r|r|r|r|r|}
\hline
\raisebox{0ex}[0.4cm][1.2\height]{Half-perimeter} & $|$Fix($1$)$|$ &
$|$Fix($r^2$)$|$ & $|$Fix($d_1$)$|$ & $|$Fix($d_2$)$|$ & \#\,Orbits &
$|$Fix($\mathfrak{D}_2$)$|$ & Asym \\
\hline \hline
 2   &   1   &   1   &   1   &   1   &   1   &   1   &   0  \\   \hline
 3   &   2   &   2   &   0   &   0   &   1   &   0   &   0  \\   \hline
 4   &   5   &   3   &   1   &   3   &   3   &   1   &   0  \\   \hline
 5   &   14   &   6   &   0   &   0   &   5   &   0   &   8   \\   \hline
 6   &   42   &   10   &   2   &   10   &   16   &   2   &   24   \\   \hline
 7   &   132   &   20   &   0   &   0   &   38   &   0   &   112   \\   \hline
 8   &   429   &   35   &   5   &   35   &   126   &   3   &   360   \\   \hline
 9   &   1430   &   70   &   0   &   0   &   375   &   0   &   1360   \\   \hline
 10   &   4862   &   126   &   14   &   126   &   1282   &   6   &   4608   \\   \hline
 11   &   16796   &   252   &   0   &   0   &   4262   &   0   &   16544   \\   \hline
 12   &   58786   &   462   &   42   &   462   &   14938   &   10   &   57840   \\   \hline
 13   &   208012   &   924   &   0   &   0   &   52234   &   0   &   207088   \\   \hline
 14   &   742900   &   1716   &   132   &   1716   &   186616   &   20   &   739376   \\   \hline
 15   &   2674440   &   3432   &   0   &   0   &   669468   &   0   &   2671008   \\   \hline
 16   &   9694845   &   6435   &   429   &   6435   &   2427036   &   35   &   9681616   \\   \hline
 17   &   35357670   &   12870   &   0   &   0   &   8842635   &   0   &   35344800   \\   \hline
 18   &   129644790   &   24310   &   1430   &   24310   &   32423710   &   70   &   129594880   \\   \hline
 19   &   477638700   &   48620   &   0   &   0   &   119421830   &   0   &   477590080   \\   \hline
 20   &   1767263190   &   92378   &   4862   &   92378   &   441863202   &   126   &   1767073824   \\   \hline
\end{tabular}
\end{footnotesize}
\caption{\small Parallelogram polyominoes enumerated by their symmetries and half-perimeter.}
\end{center}
\end{table}

\begin{table}[!ht]
\begin{center}
\begin{footnotesize}
\begin{tabular}{|c||r|r|r|r|r|r|r|}
\hline
\raisebox{0ex}[0.4cm][1.2\height]{Area} & $|$Fix($1$)$|$ & $|$Fix($r^2$)$|$ & $|$Fix($d_1$)$|$ & $|$Fix($d_2$)$|$ & \#\,Orbits &  $|$Fix($\mathfrak{D}_2$)$|$ & Asym \\
\hline \hline
1   &   1   &   1   &   1   &   1   &   1   &   1   &   0   \\   \hline
2   &   2   &   2   &   0   &   0   &   1   &   0   &   0   \\   \hline
3   &   4   &   2   &   0   &   2   &   2   &   0   &   0   \\   \hline
4   &   9   &   5   &   1   &   1   &   4   &   1   &   4   \\   \hline
5   &   20   &   4   &   0   &   4   &   7   &   0   &   12   \\   \hline
6   &   46   &   12   &   0   &   2   &   15   &   0   &   32   \\   \hline
7   &   105   &   9   &   1   &   9   &   31   &   1   &   88   \\   \hline
8   &   242   &   28   &   0   &   6   &   69   &   0   &   208   \\   \hline
9   &   557   &   21   &   1   &   21   &   150   &   1   &   516   \\   \hline
10   &   1285   &   65   &   1   &   13   &   341   &   1   &   1208   \\   \hline
11   &   2964   &   48   &   0   &   48   &   765   &   0   &   2868   \\   \hline
12   &   6842   &   152   &   2   &   32   &   1757   &   0   &   6656   \\   \hline
13   &   15793   &   111   &   1   &   111   &   4004   &   1   &   15572   \\   \hline
14   &   36463   &   351   &   1   &   73   &   9222   &   1   &   36040   \\   \hline
15   &   84187   &   257   &   3   &   257   &   21176   &   1   &   83672   \\   \hline
16   &   194388   &   814   &   2   &   172   &   48844   &   2   &   193404   \\   \hline
17   &   448847   &   593   &   3   &   593   &   112509   &   1   &   447660   \\   \hline
18   &   1036426   &   1882   &   4   &   396   &   259677   &   0   &   1034144   \\   \hline
19   &   2393208   &   1370   &   4   &   1370   &   598988   &   2   &   2390468   \\   \hline
20   &   5526198   &   4352   &   6   &   920   &   1382869   &   2   &   5520924   \\   \hline
21   &   12760671   &   3165   &   7   &   3165   &   3191752   &   1   &   12754336   \\   \hline
22   &   29466050   &   10054   &   8   &   2124   &   7369559   &   2   &   29453868   \\   \hline
23   &   68041019   &   7309   &   11   &   7309   &   17013912   &   3   &   68026396   \\   \hline
\end{tabular}
\end{footnotesize}
\caption{\small Parallelogram polyominoes enumerated by their symmetries and area.}
\end{center}
\end{table}

\end{document}